# Two-Step Many-Objective Optimal Power Flow Based on Knee Point-Driven Evolutionary Algorithm

**Yahui Li and Yang Li \***

School of Electrical Engineering, Northeast Electric Power University, Jilin 132012, China; liyh1993@gmail.com
* Correspondence: liyang@neepu.edu.cn; Tel.: +86-159-4796-6691

**Abstract:** To coordinate the economy, security and environment protection in the power system operation, a two-step many-objective optimal power flow (MaOPF) solution method is proposed. In step 1, it is the first time that knee point-driven evolutionary algorithm (KnEA) is introduced to address the MaOPF problem, and thereby the Pareto-optimal solutions can be obtained. In step 2, an integrated decision analysis technique is utilized to provide decision makers with decision supports by combining fuzzy c-means (FCM) clustering and grey relational projection (GRP) method together. In this way, the best compromise solutions (BCSs) that represent decision makers' different, even conflicting, preferences can be automatically determined from the set of Pareto-optimal solutions. The primary contribution of the proposal is the innovative application of many-objective optimization together with decision analysis for addressing MaOPF problems. Through examining the two-step method via the IEEE 118-bus system and the real-world Hebei provincial power system, it is verified that our approach is suitable for addressing the MaOPF problem of power systems.

## 1. Introduction

Optimal power flow (OPF) plays a major part role in guaranteeing the safe and economical operation of power systems [1,2], and it has been receiving the wide-spread attention of professionals and researchers from academia and industry [3,4], especially in the case of large-scale integrations of renewable energy resources [5,6]. The key idea of OPF is to find the optimal operating point with the lowest generation/operating costs under the premise of constraints [7–9], which contain a series of equality and inequality equations [10,11]. However, the conventional mono-objective OPF, which generally seeks optimum economy [12,13], such as active power losses or generation costs, becomes unable to meet the diversified needs of electricity consumers. In [12], an OPF model is proposed for determining optimal operating points of a power system, and the operating costs of the system are set to the mono-objective function in the model. In [13], the adjustable direct current OPF is presented and the objective function is taken as the total generation cost of units. And at the same time, the power flow characteristics of a modern power system are becoming increasingly complex due to the growing penetration of distributed generations [14–17] and the deployments of novel power electronic loads [18–22]. In this context, multi-objective OPF (MOPF) has received the extensive attention of researchers in the field of OPF [23–28], since it can coordinate different-weight or even conflicting multiple objectives. However, MOPF poses challenges in terms of computational complexity due to its inherent non-linear, non-convex, and non-smooth characteristic [23,24], which is hard to solve directly.

Recent research suggests that multi-objective evolutionary algorithms (MOEAs) are promising tools for addressing various challenging optimization tasks in engineering fields [29–32]. In order to optimal distributed generation planning, a MOEA is employed in [29]. Reference [30] reviews the most representative MOEAs that have been reported, and MOEA has developed as an effective method to solve such an optimization problem. In [31], MOEA is employed for planning overtime of software engineers. The layout of wind farms is optimized via MOEA in [32]. In particular, MOEAs can be also applied to solve the OPF issue [23–28]. Unfortunately, the MOPF can only cope with the optimization issue with two to three objectives, which, to a certain extent, limits the practicality of this type of methods. In addition, many-objective optimization problems (MaOPs), considering four or more objective functions in the OPF problem [33–36], are commonly existed phenomenon in the practice of real-world power system operation [33]. In [34], a specially tailored MOEA is presented for tackling the current large-scale MaOPs. Another MOEA based on adaptive search strategy is presented for coping with MaOPs in [35]. In [36], six different evolutionary algorithms (EAs) are tested, and the results prove that MOEAs exhibit their own capabilities in dealing with different MaOPs. For this reason, MaOPs have recently gained a great deal of attention as most existing MOEAs are inadequate for solving OPF problems with four or more objectives, and it has become a hotspot to enhance the ability of MOEAs for addressing MaOPs issues [37–39]. However, many-objective OPF (MaOPF) is quite challenging for solving since it is generally non-convex and non-deterministic polynomial-time hard (NP-hard). Motivated by the recent work in literature [23], a new powerful MOEA, called knee point-driven evolutionary algorithm (KnEA), is applied for solving this problem in the paper, which is helpful to better adapt the increasingly diversified operating requirements for the construction of the modern power systems.

In recent years, MOEAs have been successfully utilized in the field of multi-objective OPF (MOPF) problems in some significant pioneering works. In [24], one of MOEAs, artificial bee colony algorithm, is applied for addressing MOPF issues. For solving similar MOPF issues, the improved strength Pareto evolutionary algorithm is adopted in [25]. In [26], the gravitational search algorithm is employed to cope with this issue. In [27], MOEA is applied for solving the MOPF problems in combined heat and power economic emission dispatch. An approach based on the improved MOEA is proposed to generate Pareto-optimal solutions efficiently in [28]. In [40], a hybrid MOEA is put forward to deal with the MOPF issue by taking into account a set of various constraints. In [41], a model with two optimization objectives representing economy and stability is built for the system, then it can be solved due to the adoption of non-dominated sorting genetic algorithm II (NSGA-II) [42,43]. In [44], an improved MOEA/D algorithm is used for solving MOPF issues of power systems, and the used OPF model considers two and three objective functions in the indices relevant to cost, emissions, power losses, and stability. Unfortunately, most of all the above investigations focus on the OPF problems with two or three objectives for the power system, which is unable to meet the electricity consumers' needs of an increasingly diverse. In particular, there are many requirements that should be satisfied both for electricity suppliers and for users in the actual operation of power systems, which explains the reason why MaOPF is an urgent practical problem. As a matter of fact, the MaOPF has arisen as a consequence of some difficulty to overcome in the context of a research and development project. However, as far as the authors know, very few studies have investigated the MaOPF issue in literatures, thus this work focuses on solving the MaOPF problem.

This paper proposes a novel two-step MaOPF method by combining KnEA and integrated decision making for addressing this problem. The approach includes two-folds: many-objective optimization and decision support. At step 1, the Pareto-optimal solutions can be obtained, through solving the model of MaOPF with employing KnEA; at step 2, the best compromise solutions (BCSs) can be identified according to priority memberships in each group, and each group represents different objective preference. While the fuzzy c-means (FCM) clustering is applied to divide the solutions into several groups, grey relational projection (GRP) method is employed to calculate priority memberships. The primary contribution of the proposal is the innovative application of many-objective optimization into the optimal power flow field.

The rest of this paper is organized as follows. In Section 2, a MaOPF model is built; Section 3 displays the solution methodology based on the model; and then, Section 4 contains case studies; lastly, conclusions are drawn in Section 5.

## 2. MaOPF Model

This section outlines the model of MaOPF, including the objective functions and related constraints. With the current development of power systems, different requirements need to be met, such as economy, safe and environmental protection. Thus, four objective functions are contained in the MaOPF model. What's more, the equality and inequality constraints are also included in the novelty MaOPF model, and those common constraints in OPF issues are employed in this paper.

*2.1. Objective Function*

To satisfy the requirements of economic, safety and environmental in power system, the objective functions of MaOPF problem consider generation cost, voltage deviation, *L*-index, and emissions of polluting gases in this work.

2.1.1. Generation Costs

Generally, the minimum generation cost $f_1$ is the main objective function that must be considered in the OPF problem, which represents the economy of operation of power systems. The expression of the generation cost is [45–47]:

$$f_1 = \sum_{i=1}^{N_G} (\alpha_i P_{G,i}^2 + \beta_i P_{G,i} + \gamma_i) \quad (1)$$

where $P_{G,i}$ denotes the active output power produced by generator *i*, and $N_G$ is the total number of generators; $\alpha$, $\beta$, and $\gamma$ indicate the quadratic, linear and constant factors of a generator, respectively.

2.1.2. Index of Voltage Deviation

Taking into account that the voltage deviation is an important measure to reflect the voltage quality and safety level of a power system [24,26], the minimize voltage deviation index $f_2$ is taken as one of the optimization objectives for evaluating system security [47,48]. The expression of $f_2$ is given by:

$$f_2 = \sum_{i=1}^{N} (U_i - U_{ref,i})^2 \quad (2)$$

where $U_i$ indicates the voltage amplitude of bus *i* in the system, and the total number is *N*, $U_{ref,i}$ is the pre-defined voltage amplitude of $U_i$.

2.1.3. Static Voltage Stability Margin

For the static security issue, *L*-index is another evaluation merit in OPF problems. The value of *L*-index can judge how far from the operation point of voltage collapse to that of normal, and *L*-index is defined as [49,50]:

$$L_j = \left| 1 - \sum_{i=1}^{N_G} F_{ji} \frac{U_i}{U_j} \angle (\theta_{ij} + \delta_i - \delta_j) \right|$$
$$F_{ji} = -[\mathbf{Y}_{LL}]^{-1}[\mathbf{Y}_{LG}] \quad (3)$$

where *θ* is the phase-angle difference between two buses, and *δ* is the voltage phase angle of the bus, $\mathbf{Y}_{LL}$ and $\mathbf{Y}_{LG}$ are subarrays in the admittance matrix for nodes. Then, the objective function $f_3$, which is the minimum value of *L*-index, is employed to describe the static voltage stability margin. The equation of the third objective is shown as follows [49,50]:

$$f_3 = \max(L_j), \quad j=1,...,N_b \tag{4}$$

where $N_b$ is the number of load buses in the system.

2.1.4. Emissions of Polluting Gases

In the proposed MaOPF model, the fourth objective function $f_4$ considers the environmental demand of power system, thus the minimum emissions of polluting gases are utilized. The expression of $f_4$ is as follows [51,52]:

$$f_4 = \sum_{i=1}^{N_G}(a_i P_{G,i}^2 + b_i P_{G,i} + c_i) \tag{5}$$

where $a_i$, $b_i$ and $c_i$ denote the quadratic, linear and constant polluting gases emissions coefficients of generator *i*.

*2.2. Constraints in Power Systems*

In the power system, the used main constraints in the MaOPF model are introduced in this section.

2.2.1. Constraints of Equality

In the OPF problem, the equality constraints are universally considered and enforced, which can be written as nonlinear equations as follows [24,26]:

$$\begin{aligned}P_{g,i} - P_{d,i} &= U_i \sum_{j \in i} U_j \left(G_{ij} \sin\theta_{ij} + B_{ij} \cos\theta_{ij}\right) \\ Q_{g,i} - Q_{d,i} &= U_i \sum_{j \in i} U_j \left(G_{ij} \sin\theta_{ij} - B_{ij} \cos\theta_{ij}\right)\end{aligned} \tag{6}$$

where $P_{g,i}$ and $P_{d,i}$ are the injected active power and active load in bus *i*, while $Q_{g,i}$ and $Q_{d,i}$ are the reactive injected power and reactive load, $G_{ij}$ and $B_{ij}$ are respectively the conductance and susceptance between buses *i* and *j*. This equation suggests that the active and reactive powers need to keep balance in the power system.

2.2.2. Constraints of Inequality

Herein, the bounds of variables in the power system are considered for the purpose of ensuring the power system in a safe state during operation, and they can be formulated as [24,26]:

$$\begin{aligned}P_{G,i}^{\min} &\leq P_{G,i} \leq P_{G,i}^{\max}, \quad i=1,...,N_G \\ Q_{G,i}^{\min} &\leq Q_{G,i} \leq Q_{G,i}^{\max}, \quad i=1,...,N_G \\ U_i^{\min} &\leq U_i \leq U_i^{\max}, \quad i=1,...,N \\ T_i^{\min} &\leq T_i \leq T_i^{\max}, \quad i=1,...,N_T \\ Q_{C,i}^{\min} &\leq Q_{C,i} \leq Q_{C,i}^{\max}, \quad i=1,...,N_C \\ S_{L,i}^{\min} &\leq S_{L,i} \leq S_{L,i}^{\max}, \quad i=1,...,N_L\end{aligned} \tag{7}$$

where $Q_{G,i}$ represents the reactive power produced by generator *i*; $T_i$ indicates the tap of adjustable transformer *i*, and the number of transformers is $N_T$; $Q_{C,i}$ expresses the reactive power of compensation device *i*, and the number of devices is $N_C$; $S_{L,i}$ denotes the power flow in the branch $i$ of the system, and $N_L$ expresses the number of branches; each variable has its upper and lower limits, which are respectively represented by the subscript 'max' and 'min'.

## 3. Two-Step Solution Approach

The proposal is divided into two phases: an optimization process with many objectives and the following decision support procedure. At step 1, the set of Pareto optimal solutions is gained with the employment of KnEA. Then at step 2, FCM clustering is adopted to classify the set obtained in the first step. GRP method is applied to automatically select BCSs from each group.

*3.1. Optimization Process with Many Objectives*

The first step is introduced in this section. By solving MaOPF model with KnEA, the set of Pareto optimals is obtained.

3.1.1. KnEA-Based Many-Objective Optimization

The key concept of the KnEA is knee points (KPs). Different from other MOEAs, the KnEA employs a prevailing non-dominance selection criterion strategy and a secondary selection criterion strategy whose reference is the KPs in the optimization process [33]. After the environmental selection process, the diversity of the population is improved. In addition, mating selection adopting tournament strategy is applied in KnEA, and weighted distance is a metric of the strategy. The weighted distance $WD(p)$ of solution p is defined as [33]:

$$WD(p) = \sum_{i=1}^{k} wd_{p_i} dis_{pp_i}$$
$$s.t. \quad wd_{p_i} = \frac{rd_{p_i}}{\sum_{i=1}^{k} rd_{p_i}} \quad (8)$$
$$rd_{p_i} = \frac{1}{\left| dis_{pp_i} - \frac{1}{k}\sum_{i=1}^{k} dis_{pp_i} \right|}$$

where $p_i$ is the *i*th nearest neighbor of solution $p$; $wd_{p_i}$ and $rd_{p_i}$ are respectively the weight and rank of $p_i$; $dis_{pp_i}$ denotes the Euclidean distance between $p$ and $p_i$.

In the KnEA, a KP is defined as the point with the maximum distance to the hyperplane within neighborhood scope in the objective functions space [53]. In the *g*th generation, the neighborhood scope $R_g^i$ corresponding to the objective function $f_i$ is given by:

$$R_g^i = \left( f_{i,g}^{\max} - f_{i,g}^{\min} \right) \cdot r_g$$
$$s.t. \quad r_g = r_{g-1} \cdot e^{-\frac{1 - t_{g-1}/TH}{N_{obj}}} \quad (9)$$

where the upper and lower limits of the *i*th objective function $f_i$ are respectively represented by the superscript 'max' and 'min', and the total number of functions is $N_{obj}$; $r$ denotes the proportion of neighborhood scope in the objective span; $g$ represents the *g*th generation; $t$ is the proportion of the KPs in the whole population; $0<TH<1$ expresses the boundary.

3.1.2. Procedure of Many-Objective Optimization

The flowchart of the optimization procedure using the KnEA is shown in Figure 1, and the specific steps are as follows [33].

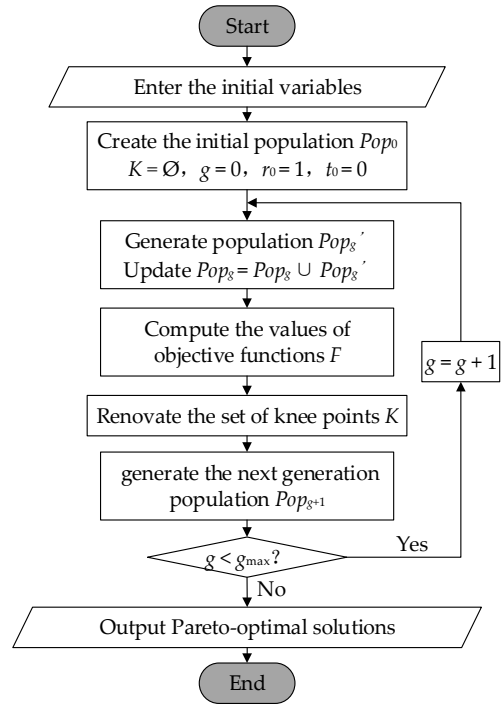

**Figure 1.** Process of optimization scheme by using knee point-driven evolutionary algorithm (KnEA).

Step 1: Enter the initial variables. The variables mainly incorporate three-folds as follows.

(1) The network parameters: the related information of power systems.
(2) The controlled variable parameters: the bounds which are shown as (7), and the steps of $T$ and $Q_C$. The considering controlled variable are listed as follow:

$$\left( \underbrace{P_{G,1}, \cdots, P_{G,i}, \cdots, P_{G,N_G}, U_{G,1}, \cdots, U_{G,i}, \cdots, U_{G,N_G}}_{\text{continuous variables}}, \underbrace{T_1, \cdots, T_i, \cdots, T_{N_T}, Q_{C,1}, \cdots, Q_{C,i}, \cdots, Q_{C,N_C}}_{\text{discrete variables}} \right)$$

where $U_G$ is the generator terminal voltage.

(3) The algorithm parameters: the population size $N_{pop}$, the maximum generation number $g_{max}$, the set of KPs $K$, the ratio of size $r$, the rate of KPs in population $t$, the number of objectives $N_{obj}$ which is taken as 4 in this paper.

Step 2: Create an initial population $Pop_0$, initial the set of KPs $K$ which is an empty set, and set the generation counter $g=0$, $r_0=1$, $t_0=0$ Considering that the variables $T$ and $Q_C$ are discrete variables while others are continuous variables, a hybrid coding scheme is employed in the initialization process.

Step 3: Generate population $Pop_g'$ by adopting binary tournament mating selection with three strategies for distinguishing solutions, that is, dominance comparison, KP criterion, and $WD(p)$ in Equation (8) [33]. If the non-dominance selection and the secondary selection criterion strategy fail to discriminate solutions in $Pop_g'$, they can be chosen eventually according to the value of $WD(p)$. Then, the population $Pop_g'$ is formed by the solutions with the mating selection.

Step 4: Update the population $Pop_g$ with genetic variations. The two operations, simulated binary crossover, and polynomial mutation, are implemented in $Pop_g'$. The population $Pop_g$ is then updated based on the individuals in $Pop_g$ and $Pop_g'$, and $Pop_g = Pop_g \cup Pop_g'$.

Step 5: Compute the values of objective functions $F = \{f_1, f_2, f_3, f_4\}$, and then the non-dominated solutions are identified from $Pop_g$.

Step 6: Renovate the set $K$. The KPs are chosen according to (9), and the set of KPs $K$ are recorded and renovated in the optimization process.

Step 7: Based on the set $K$ and the objective function values, generate the next generation population $Pop_{g+1}$ based on the environmental selection strategy [33].

Step 8: Judge the termination criteria. If $g < g_{max}$, return to Step 3 after adding 1 to the current generation *g*; otherwise, output the Pareto-optimal solutions.

*3.2. Decision Support*

In real-world practice, it is quite challenging for decision makers to figure out whether a Pareto-optimal solution is a BCS or not from among the non-inferior solutions. First, considering the fact that there are a lot of generated Pareto-optimal solutions, the references of decision makers might be different for a specific operation point. Another issue is that for a specific system the preference of the same decision maker may also vary according to changing operational requirements. Therefore, in the second step, that is, the decision support process, FCM and GRP method are adopted to evaluating Pareto optimal solutions, and BCSs are identified which are represent decision-makers' different, even conflicting, preferences.

3.2.1. Fuzzy *c*-Means

FCM is a classical unsupervised clustering algorithm. After processing of FCM clustering, the similarities between the solutions in the same group are the largest, while the similarities between the solutions in different groups are the smallest. The model of FCM can be formulated as [54–56]:

$$\min \ J_m = \sum_{i=1}^{N_p} \sum_{j=1}^{N_c} \mu_{ij}^m \|w_i - c_j\|^2$$
$$s.t. \ \sum_{j=1}^{N_c} \mu_{ij} = 1 \tag{10}$$

where $J$ expresses loss function for judging the convergent degree, $m(m \geq 1)$ is a constant for controlling the clustering fuzziness, $\mu_{ij}(\mu_{ij} \in [0,1])$ is the membership degree between solution $w_i$ and center $c_j$, while $w_i$ denotes the *i*th Pareto optimal solution in the whole set, $c_j$ is the *j*th clustering center; $N_p$ and $N_c$ are the number of solutions and clusters, respectively. Here, $N_c$ is taken as 4 (corresponding to the considered four objective functions).

3.2.2. GRP Method

As one of the decision methods, GRP method is especially suitable for evaluating the solutions with grey relationships in the MaOPs, which is based on grey system theory and vector projection [23]. The projection $V_l^{+(-)}$ of the *l*th solution in the ideal scheme can be written as [57,58]:

$$V_l^{+(-)} = \sum_{k=1}^{N_p} \gamma_{lk}^{+(-)} \frac{\omega_k^2}{\sqrt{\sum_{k=1}^{N_l}(\omega_k)^2}} \tag{11}$$

where plus sign denotes positive scheme, minus sign represents negative scheme, $\gamma_{lk}$ indicates the grey relational factor of the $k$th objective in $l$th solution. $\omega_k$ is the weight of the $k$th objective, the corresponding weights of four objectives are set to the same value in this paper, and the operators can adjust the weights according to the actual working condition or personal preference. Then, the priority membership $PM_l\,(0 \leq PM_l \leq 1)$ of solution $l$ can be written as follows:

$$PM_l = \frac{(V_0 - V_l^-)^2}{(V_0 - V_l^-)^2 + (V_0 - V_l^+)^2} \tag{12}$$

where $V_0$ equals to the value of $V_l$ when $V_l$ takes 1. The greater the membership of the solutions is, the closer it is to the ideal scheme; and vice versa. In this way, the solutions with the highest PM values in each group are regarded as the BCSs.

## 4. Case Studies

For examining validity of the approach provided in this paper, two test systems with varied complexity levels, i.e., the IEEE standard system and the system applied in Hebei province, are taken as test cases. And furthermore, to properly measure the optimization performance of our approach, two state-of-the-art MOEAs for solving MaOPs, i.e., the reference vector guided evolutionary algorithm (RVEA) and non-dominated sorting genetic algorithm III (NSGA-III), are employed as comparison algorithms. All programs in this work are carried out by a desktop computer with 3.40 GHz CPU basic frequency and 4 GB memory.

### 4.1. IEEE 118-Bus System

The first step is introduced in this section. By solving MaOPF model with KnEA, the set of Pareto optimals is obtained.

#### 4.1.1. Introduction to the System

As a well-known test system, IEEE 118-bus system is extensively studied in previous literature [26]. This system with base capacity 100 MVA includes 14 active generators, 132 branches, 9 adjustable transformers. In the system, bus 69 is the slack bus. The related coefficients of generator $i$, such as $\alpha_i$, $\beta_i$, $\gamma_i$, $a_i$, $b_i$, and $c_i$, are extracted from literature [26].

The limits of controlled variables are listed as follows: the lower and upper bounds of the voltage are respectively 0.95 p.u. and 1.10 p.u., the tap $T$ varies from 0.9 p.u. to 1.1 p.u., and the lower and upper bounds of $Q_C$ are 0 and 0.5 p.u.; the step-size of $T$ and $Q_C$ are respectively 0.0125 p.u. and 0.01 p.u.; the upper bound of the branch transmission capacity is 300 MVA.

#### 4.1.2. Algorithm Comparison

As mentioned above, three algorithms are employed to solve MaOPF problem. In order to facilitate comparison, $N_{pop}$ and $g_{max}$ of all three algorithms are respectively 50 and 100. The three algorithms repeatedly run 20 times independently. Among all the results of the 20 runs for each algorithm, without loss of generality, one result (i.e., a set of Pareto optimals) is randomly taken as an instance for the consequent analysis. The distributions of three selected results with four objective functions are shown in Figures 2–4.

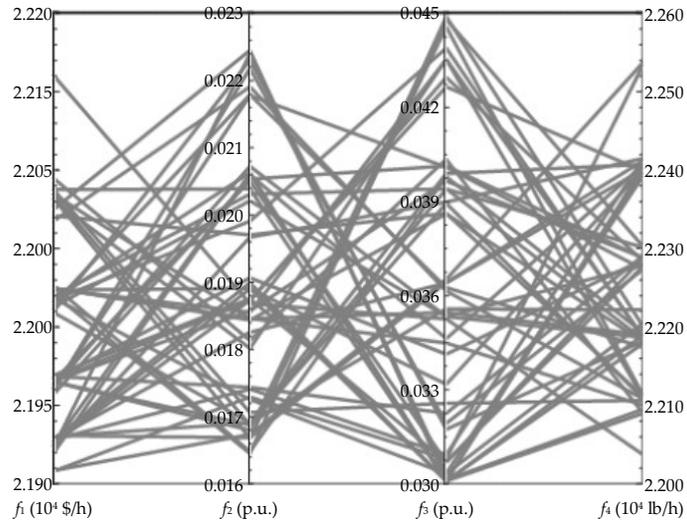

**Figure 2.** Distribute condition of Pareto-optimal solutions of KnEA.

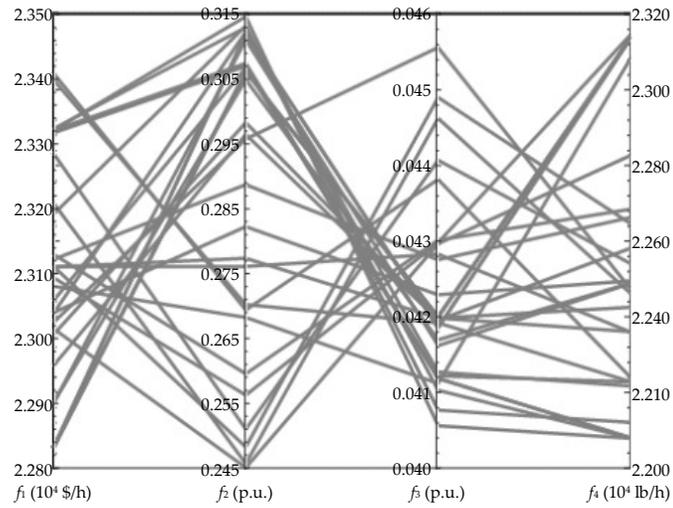

**Figure 3.** Distribute condition of Pareto-optimal solutions of reference vector guided evolutionary algorithm (RVEA).

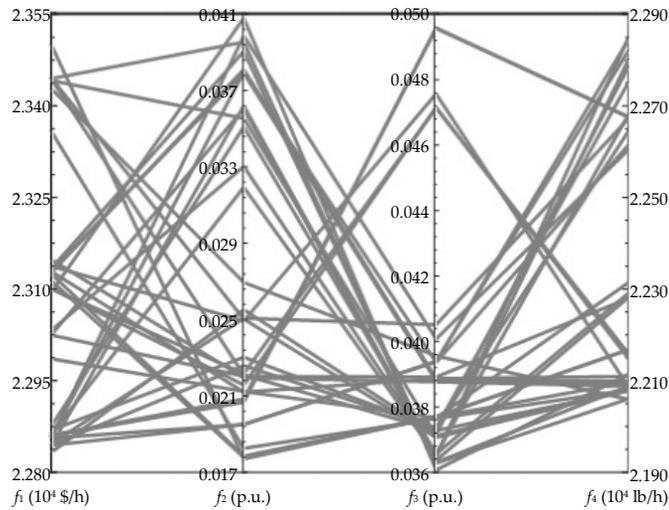

**Figure 4.** Distribution of Pareto-optimal solutions of non-dominated sorting genetic algorithm III (NSGA-III).

The comparison of extreme values of the four objective functions obtained by each algorithm is shown in Table 1, and the smallest values in each line of Table 1 have been marked with bold fonts. According to the results in Table 1, it is obvious that the extreme values of the four objectives obtained by KnEA are smaller than that of the other two algorithms when solving the MaOPF problem. To a certain extent, the extreme value can evaluate the performance of MOEAs. The smaller extreme value means that the optimization performance is more effective, thus KnEA has the best performance in the three algorithms only from the view of extreme values.

**Table 1.** Comparison of extreme values of KnEA, RVEA, and NSGA-III. (KnEA: knee point-driven evolutionary algorithm; RVEA: reference vector guided evolutionary algorithm; NSGA-III: non-dominated sorting genetic algorithm III).

| Objective Function | Extreme Value | KnEA | RVEA | NSGA-III |
|---|---|---|---|---|
| $f_1/(10^4$ \$/h$)$ | Maximum value | **2.3062** | 2.3413 | 2.3507 |
|  | Minimum value | **2.2808** | 2.2826 | 2.2831 |
| $f_2/$(p.u.) | Maximum value | **0.0224** | 0.3147 | 0.0408 |
|  | Minimum value | **0.0165** | 0.2451 | 0.0177 |
| $f_3/$(p.u.) | Maximum value | **0.0449** | 0.0456 | 0.0496 |
|  | Minimum value | **0.0301** | 0.0406 | 0.0361 |
| $f_4/(10^4$ lb/h$)$ | Maximum value | **2.2539** | 2.3151 | 2.2857 |
|  | Minimum value | **2.2036** | 2.2078 | 2.2058 |

How to assess the performance of MOEAs has recently been attracting concerns. Unfortunately, this is still an open question at the moment. In general, a good evaluation indicator should have good convergence and distribution characteristic [23]. Two quantitative indicators, which can assess the optimization performances of three different algorithms in different aspects, are employed in this study.

(1) Generational distance

The first indicator is the well-known generational distance (GD), which represents the convergence conditions of the set [23]. For measuring the convergence of obtained solutions, the formulation of GD is given as follows:

$$\text{GD} = \frac{\sqrt{\sum_{i=1}^{N_p} D_i^2}}{N_p} \quad (13)$$

where $D_i$ denotes the Euclidean distance in objective function space, which is calculated between each two nearest solutions.

(2) Spacing

The spacing (SP) is another popular indication for estimating the distribution of a Pareto front, and its expression is given by [23]:

$$\text{SP} = \sqrt{\frac{1}{N_p - 1} \sum_{i=1}^{N_p} (\bar{D} - D_i)^2} \quad (14)$$

where $\bar{D}$ represents the average value of $D_i$. It should be noted that a solution with smaller values of the above two metrics has better performances about convergence and diversity.

In view of the randomness of MOEAs to optimal results [23], all the used three algorithms are independently carried out 20 times. In Table 2, the obtained best, average and worst values of two metrics are listed.

**Table 2.** Statistical values of two metrics for the three algorithms (GD: generational distance; SP: spacing).

| Algorithm | Metrics | Best | Average | Worst |
|---|---|---|---|---|
| KnEA | GD | 4133.68 | 4515.35 | 4868.10 |
| | SP | 15.23 | 16.40 | 17.92 |
| RVEA | GD | 5347.71 | 5893.61 | 6286.75 |
| | SP | 40.16 | 65.99 | 69.43 |
| NSGA-III | GD | 4879.06 | 5430.93 | 6250.68 |
| | SP | 17.37 | 19.67 | 21.35 |

From Table 2, we can see that the metrics GD and SP of the KnEA are better than those of the RVEA and NSGA-III. These results suggest that the KnEA has advantages over the other alternatives in the convergence and distribution performances.

Subsequently, the average calculation times of the three algorithms in 20 runs are presented in Table 3. And it is easy to judge from the average times in Table 3, comparing with RVEA and NSGA-III, the optimization speed of KnEA for solving MaOPF problem is better.

**Table 3.** Average times of each algorithm.

| Algorithm | KnEA | RVEA | NSGA-III |
|---|---|---|---|
| Average Time (s) | 88.12 | 93.45 | 90.69 |

From the above comparison, it is clear that KnEA is superior to RVEA and NSGA-III in optimization effects and solution efficiency in solving the MaOPF problem.

4.1.3. Result Analysis

Taking the representative Pareto-optimal solutions obtained by KnEA as an example, the solutions are divided into four groups, which corresponds to the four objective functions, through the FCM clustering, and the distribution of four groups of solutions is shown via different colors in Figure 5.

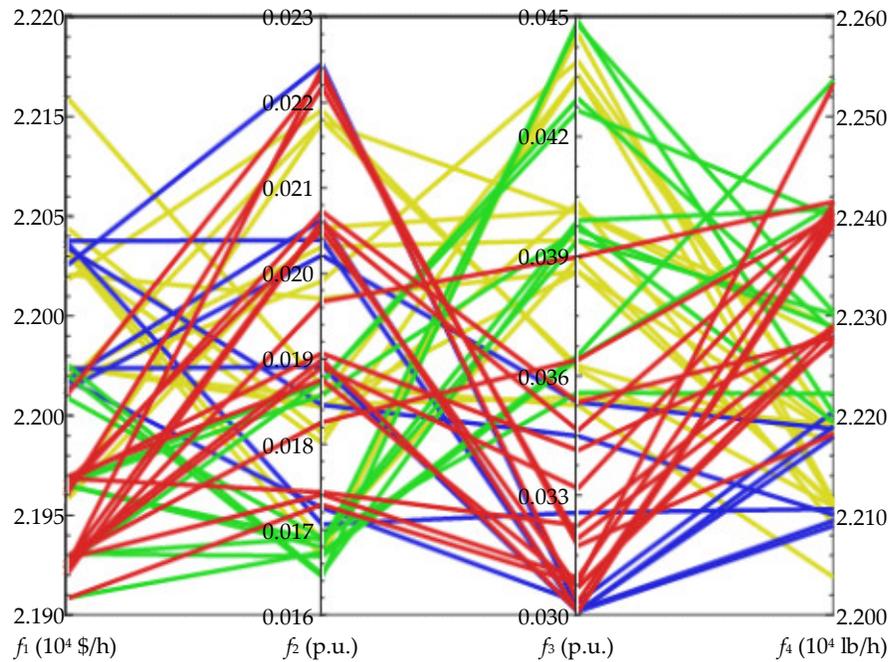

**Figure 5.** Distribution of Pareto-optimal solutions of KnEA after Fuzzy c-Means (FCM) clustering.

In Figure 5, each line denotes one solution in the set, and the lines with red, green, blue and yellow colors represent that decision makers prefer for $f_1$, $f_2$, $f_3$, $f_4$, respectively. When many lines cross between two adjacent objectives, it indicates that the two objectives are in a conflicting relationship. It should be noted that each Pareto-optimal solution acquired from KnEA is not the best for every objective since, for a MaOPF problem, they are only non-inferior solutions.

GRP method is used to evaluate the solutions after adopting FCM clustering, and each group belongs to one scheme. After the membership of each solution is computed, the BCSs, which have the highest membership values in each group, are listed in Table 4.

**Table 4.** Best compromise solutions (BCSs) of IEEE 118-bus system.

| BCSs | $f_1$ ($10^4$ \$/h) | $f_2$ (p.u.) | $f_3$ (p.u.) | $f_4$ ($10^4$ lb/h) | PM |
|---|---|---|---|---|---|
| Prefer for $f_1$ | 2.2828 | 0.0174 | 0.0309 | 2.2417 | 0.7553 |
| Prefer for $f_2$ | 2.2831 | 0.0167 | 0.0364 | 2.2413 | 0.7423 |
| Prefer for $f_3$ | 2.2919 | 0.0173 | 0.0303 | 2.2206 | 0.6885 |
| Prefer for $f_4$ | 2.2926 | 0.0185 | 0.0355 | 2.2190 | 0.6848 |

According to the BCSs shown in Table 4, the two-step mean is capable of addressing the MaOPF problem. Not only a complete and evenly distributed set is achieved, but also the BCSs can be identified.

The BCS prefer for $f_1$ is an example, and comparison results of generator variables are displayed in Table 5, before and after adopting KnEA. Furthely, the comparison of objective functions before optimization and BCS prefer for $f_1$ are shown in Table 6.

**Table 5.** Comparison results of generator variables.

| Generators | Before Optimization | | | After Optimization | | |
|---|---|---|---|---|---|---|
| | $P_G$ (p.u.) | $Q_G$ (p.u.) | $U_G$ (p.u.) | $P_G$ (p.u.) | $Q_G$ (p.u.) | $U_G$ (p.u.) |
| G1 | 4.500 | 0 | 1.050 | 4.471 | −0.856 | 1.019 |
| G2 | 0.850 | 0 | 0.990 | 0.935 | 0.489 | 0.987 |
| G3 | 2.200 | 0 | 1.050 | 2.420 | 1.811 | 1.015 |
| G4 | 3.140 | 0 | 1.015 | 3.454 | −1.897 | 1.004 |
| G5 | 2.040 | 0 | 1.025 | 2.244 | 0.331 | 1.009 |
| G6 | 0.480 | 0 | 0.955 | 0.528 | 0.301 | 0.983 |
| G7 | 1.550 | 0 | 0.985 | 1.705 | 0.793 | 1.005 |
| G8 | 1.600 | 0 | 0.995 | 1.760 | −0.310 | 1.002 |
| G9 | 3.910 | 0 | 1.005 | 4.300 | 3.749 | 1.000 |
| G10 | 3.920 | 0 | 1.050 | 4.312 | −3.919 | 1.019 |
| G11 | 5.164 | 0 | 1.035 | 3.690 | −0.772 | 1.031 |
| G12 | 4.770 | 0 | 1.040 | 4.501 | 0.026 | 1.019 |
| G13 | 6.070 | 0 | 1.005 | 5.463 | −0.154 | 1.015 |
| G14 | 2.520 | 0 | 1.017 | 2.772 | 0.304 | 1.008 |

**Table 6.** Comparison results before and after optimization.

| Optimization Condition | $f_1$ ($10^4$ $/h) | $f_2$ (p.u.) | $f_3$ (p.u.) | $f_4$ ($10^4$ lb/h) |
|---|---|---|---|---|
| Before optimization | 13.1221 | 0.0416 | 1.8729 | 2.9153 |
| After optimization | 2.2828 | 0.0174 | 0.0309 | 2.2417 |

From the above table, it can be seen that the variables are all in the predefined range, and the distribution of power flow becomes more reasonable through optimization, which embodies the four objective functions after optimization are superior to their corresponding values before optimization. Accordingly, it can be concluded that the presented algorithm is an efficient tool to determine the BCSs for the MaOPF problem, which helps to provide more realistic options representing decision makers' different references.

*4.2. Application to the Hebei Provincial System*

The two-step mean is employed to an actual physical power system to evaluate the applicability, and the Hebei provincial power system located in China is further tested in this paper. This system contains 45 active generators, 169 substations with voltage grades 220 kV and above, and some compensation equipment [29]. In addition, the system has 17 channels which can extended to other power systems.

Given the maximum generation number $g_{max} = 150$, the optimization results via the approach are shown in Figure 6.

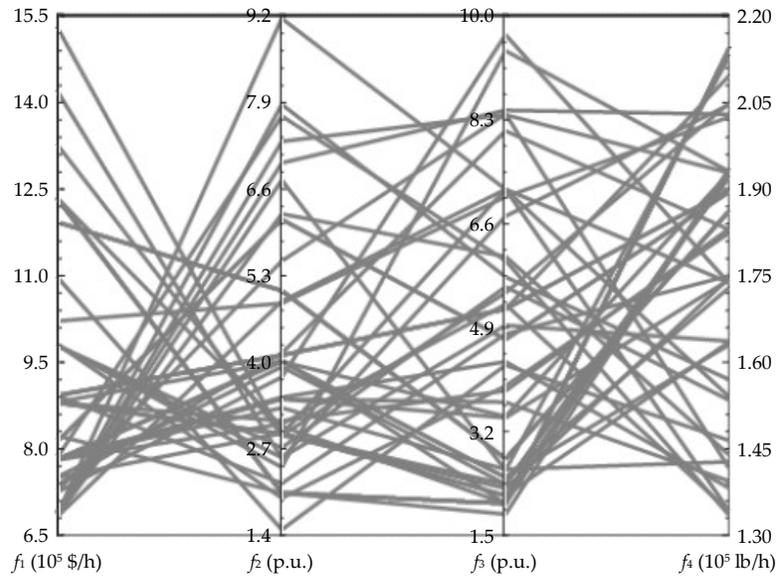

**Figure 6.** Distribute condition of solutions of the Hebei provincial system.

And then, FCM is applied for clustering the solutions, which is acquired by KnEA, into four groups, and the distributions using different colors is illustrated in Figure 7.

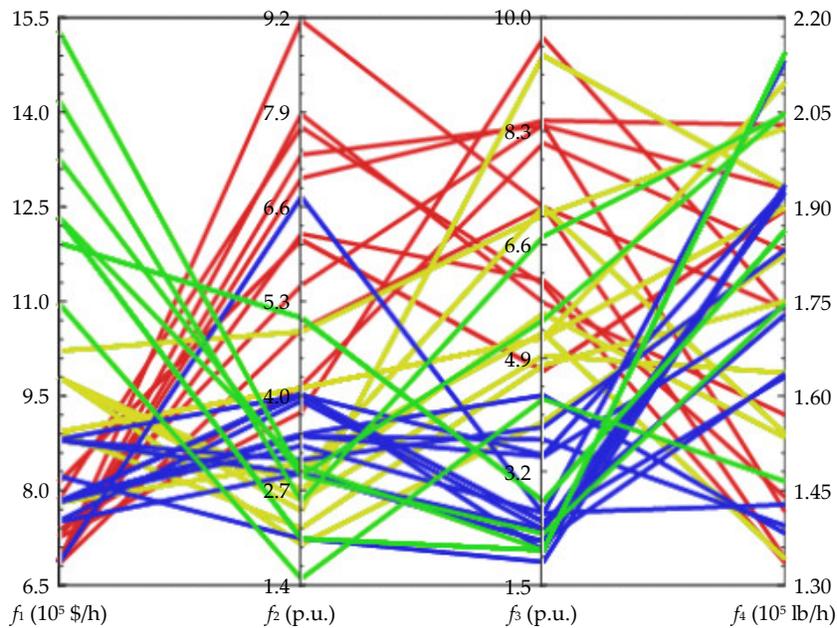

**Figure 7.** Distribute condition of Solutions of Hebei provincial system after FCM clustering.

Similar to Figure 5, the lines with different colors in Figure 7 denote different preferences of decision makers. Table 7 shows the BCSs with the maximum memberships in four groups after using GRP method.

**Table 7.** BCSs of the Hebei provincial power system.

| BCSs | $f_1$ (10$^5$ $/h) | $f_2$ (p.u.) | $f_3$ (p.u.) | $f_4$ (10$^5$ lb/h) | PM |
|---|---|---|---|---|---|
| Prefer for $f_1$ | 6.8019 | 4.9166 | 7.1676 | 1.3378 | 0.8082 |
| Prefer for $f_2$ | 11.0076 | 1.4871 | 4.2874 | 1.4632 | 0.8471 |
| Prefer for $f_3$ | 7.8233 | 3.3164 | 2.5808 | 1.4277 | 0.7163 |
| Prefer for $f_4$ | 7.8030 | 3.0755 | 5.2818 | 1.3366 | 0.8213 |

For purpose of assessing the optimization effects of our approach, the result before optimization and the obtained BCS preferring to $f_1$ are listed in Table 8.

**Table 8.** Comparison results before and after optimization.

| Optimization condition | $f_1$ (10$^5$ $/h) | $f_2$ (p.u.) | $f_3$ (p.u.) | $f_4$ (10$^5$ lb/h) |
|---|---|---|---|---|
| Before optimization | 16.9414 | 166.2188 | 70.7107 | 6.5995 |
| After optimization | 6.8019 | 4.9166 | 7.1676 | 1.3378 |

From the above table, it is clear that all the four objective functions have been improved through the proposed KnEA-based two-step MaOPF approach. Therefore, the conclusion can be drawn safely that the two-step approach is also suitable for addressing the MaOPF problems in a real-world power system.

*4.3. Discussions*

From the results, the MaOPF issues with more than three objective functions can be effectively solved both in the IEEE standard power systems and in the actual power systems. Meanwhile, KnEA is selected as the most effective algorithm in three MOEAs by comparing evaluation indicators of optimization performance. However, there are still some limitations of the performed work. As an important basic theory research, this work meant to solve the MaOPF issues in practical power systems. A simple OPF model is employed in this paper, and traditional constraints are used. A more practical OPF model will be explored in the future to consider real-world demands, such as dynamic security [41,59], and the reactive power and voltage magnitude [60]. Aimed at security problems in the power system, two safety-related functions are contained in the OPF model. $N-1$ security constraints of power systems need also to be considered for preventive and corrective actions. What's more, more static and dynamic security functions and constraints can be added to the MaOPF model for ensuring the safe and stable operation of the power system. Moreover, the configuration of static var compensation devices is also a practical problem in planning, designing, and operation, and this issue is of great significance.

**5. Conclusions**

A two-step MaOPF approach using KnEA algorithm is presented in this paper. According to the analysis of the IEEE 118-bus system and a real-world power system (i.e., Hebei provincial system, China), the following conclusions are safely drawn:

(1) Considering the generation cost, voltage deviation, static voltage stability margin and emissions of polluting gases, a MaOPF model is proposed to better adapt the increasingly diversified operating requirements of power systems.
(2) The proposed solution approach not only can yield multiple well-distributed set of Pareto-optimal solutions, but also can further determine BCSs from each group, which represent decision-makers' different, even conflicting, preferences.
(3) The simulation results on two test cases with varied complexity levels verify the effectiveness of the proposal. More importantly, the KnEA has significant advantages in the optimization performance, compared with the other popular algorithms, such as RVEA and NSGA-III.

In our future research, distributed and parallel computing techniques will be employed to further improve the problem-solving efficiency of the proposed approach. Furthermore, the definition and validation of performance metrics for multi-objective evolutionary algorithms is an unsolved very important issue nowadays. The OPF with energy storage is another beneficial topic for future study [61,62].

**Author Contributions:** Conceptualization, Y.L. (Yang Li); Methodology, Y.L. (Yang Li) and Y.L. (Yahui Li); Software, Y.L. (Yahui Li); Validation, Y.L. (Yang Li) and Y.L. (Yahui Li); Formal Analysis, Y.L. (Yang Li) and Y.L. (Yahui Li); Investigation, Y.L. (Yang Li) and Y.L. (Yahui Li); Resources, Y.L. (Yang Li); Data Curation, Y.L. (Yang Li); Writing—Original Draft Preparation, Y.L. (Yahui Li); Writing—Review and Editing, Y.L. (Yang Li); Visualization, Y.L. (Yang Li) and Y.L. (Yahui Li); Supervision, Y.L. (Yang Li); Project Administration, Y.L. (Yang Li); Funding Acquisition, Y.L. (Yang Li).

**Funding:** This research was funded by the China Scholarship Council (CSC) under Grant No. 201608220144.

**Conflicts of Interest:** The authors declare no conflict of interest.

**Abbreviations**

| | |
|---|---|
| OPF | Optimal power flow |
| MOPF | Multi-objective optimal power flow |
| MaOPF | Many-objective optimal power flow |
| MOEA | Multi-objective evolutionary algorithm |
| MaOP | Many-objective optimization problem |
| NSGA-III | Non-dominated sorting genetic algorithm III |
| KnEA | Knee point-driven evolutionary algorithm |
| RVEA | Reference vector guided evolutionary algorithm |
| BCS | Best compromise solution |
| FCM | Fuzzy c-means |
| GRP | Grey relational projection |
| $P_G$ | Active power output of a generator |
| $Q_G$ | Reactive power output of a generator |
| $N_G$ | The number of generators |
| $U$ | Voltage amplitude of a bus |
| $U_{ref}$ | Reference voltage amplitude of a bus |
| $N$ | The number of buses |
| $\theta$ | Phase-angle difference between two buses |
| $\delta$ | Voltage phase angle of a bus |
| $N_b$ | The number of load buses |
| $P_g$ | Injected active power of a load bus |
| $Q_g$ | Injected reactive power of a load bus |
| $P_d$ | Active loads of a load bus |
| $Q_d$ | Reactive load sof a load bus |
| $T$ | The tap of a transformer |
| $N_T$ | The number of adjustable transformer taps |
| $Q_C$ | The switching capacity of a reactive power compensation capacitor |
| $N_C$ | The number of reactive power compensation capacitors |
| $S_L$ | The power flow in the branch |
| $N_L$ | The number of branches |
| $WD$ | The weighted distance of solutions |
| $PM$ | The priority membership of solutions |
| $GD$ | The generational distance |
| $SP$ | The spacing |

**References**

1. Dommel, H.W.; Tinney, W.F. Optimal power flow solutions. *IEEE Trans. Power Appl. Syst.* **1968**, *10*, 1866–1876, doi:10.1109/tpas.1968.292150.
2. Vaccaro, A.; Canizares, C.A. An affine arithmetic-based framework for uncertain power flow and optimal power flow studies. *IEEE Trans. Power Syst.* **2017**, *32*, 274–288, doi:10.1109/tpwrs.2016.2565563.
3. Yang, Z.; Zhong, H.; Bose, A.; Xia, Q.; Kang, C. Optimal power flow in AC–DC grids with discrete control devices. *IEEE Trans. Power Syst.* **2018**, *33*, 1461–1472, doi:10.1109/tpwrs.2017.2721971.
4. Yang, Z.; Zhong, H.; Xia, Q.; Kang, C. A novel network model for optimal power flow with reactive power and network losses. *Electr. Power Syst. Res.* **2017**, *144*, 63–71, doi:10.1016/j.epsr.2016.11.009.
5. Li, Y.; Yang, Z.; Li, G.; Zhao, D.; Tian, W. Optimal scheduling of an isolated microgrid with battery storage considering load and renewable generation uncertainties. *IEEE Trans. Ind. Electron.* **2019**, *66*, 1565–1575, doi:10.1109/tie.2018.2840498.
6. Li, Y.; Yang, Z.; Li, G.; Mu, Y.; Zhao, D.; Chen, C.; Shen, B. Optimal scheduling of isolated microgrid with an electric vehicle battery swapping station in multi-stakeholder scenarios: A bi-level programming approach via real-time pricing. *Appl. Energy* **2018**, *232*, 54–68, doi:10.1016/j.apenergy.2018.09.211.
7. Yang, Z.; Zhong, H.; Xia, Q.; Bose, A.; Kang, C. Optimal power flow based on successive linear approximation of power flow equations. *IET Gener. Transm. Dis.* **2016**, *10*, 3654–3662, doi:10.1049/iet-gtd.2016.0547.
8. Bagde, B.Y.; Umre, B.S.; Dhenuvakonda, K.R. An efficient transient stability-constrained optimal power flow using biogeography-based algorithm. *Int. Trans. Electr. Energy Syst.* **2018**, *28*, e2467, doi:10.1002/etep.2467.
9. Sun, J.; Li, Y. Social cognitive optimization with tent map for combined heat and power economic dispatch. *Int. Trans. Electr. Energy Syst.* **2018**, e2660, doi:10.1002/etep.2660.
10. Zhang, Y.; Shen, S.; Mathieu, J.L. Distributionally robust chance-constrained optimal power flow with uncertain renewables and uncertain reserves provided by loads. *IEEE Trans. Power Syst.* **2017**, *32*, 1378–1388, doi:10.1109/tpwrs.2016.2572104.
11. Dall'Anese, E.; Baker, K.; Summers, T. Chance-constrained AC optimal power flow for distribution systems with renewables. *IEEE Trans. Power Syst.* **2017**, *32*, 3427–3438, doi:10.1109/tpwrs.2017.2656080.
12. Sharifzadeh, H.; Amjady, N. Stochastic security-constrained optimal power flow incorporating preventive and corrective actions. *Int. Trans. Electr. Energy Syst.* **2016**, *26*, 2337–2352, doi:10.1002/etep.2207.
13. Wen, Y.; Guo, C. Adjustable risk-based direct current optimal power flow. *Int. Trans. Electr. Energy Syst.* **2015**, *25*, 3212–3226, doi:10.1002/etep.2031.
14. Wang, Y.; Zhang, N.; Chen, Q.; Yang, J.; Kang, C.; Huang, J. Dependent discrete convolution based probabilistic load flow for the active distribution system. *IEEE Trans. Sustain. Energy* **2017**, *8*, 1000–1009, doi:10.1109/tste.2016.2640340.
15. Lu, Z.; Li, H.; Qiao, Y. Probabilistic flexibility evaluation for power system planning considering its association with renewable power curtailment. *IEEE Trans. Power Syst.* **2018**, *33*, 3285–3295, doi:10.1109/tpwrs.2018.2810091.
16. Jin, P.; Li, Y.; Li, G.; Chen, Z.; Zhai, X. Optimized hierarchical power oscillations control for distributed generation under unbalanced conditions. *Appl. Energy* **2017**, *194*, 343–352, doi:10.1016/j.apenergy.2016.06.075.
17. Peng, Q.; Low, S.H. Distributed optimal power flow algorithm for radial networks, I: Balanced single phase case. *IEEE Trans. Smart Grid* **2018**, *9*, 111–121, doi:10.1109/tsg.2016.2546305.
18. Xia, M.; Lai, Q.; Zhong, Y.; Li, C.; Chiang, H.D. Aggregator-based interactive charging management system for electric vehicle charging. *Energies* **2016**, *9*, 159, doi:org/10.3390/en9030159.
19. He, H.; Guo, J.; Peng, J.; Tan, H.; Sun, C. Real-time global driving cycle construction and the application to economy driving pro system in plug-in hybrid electric vehicles. *Energy* **2018**, *152*, 95–107, doi:10.1016/j.energy.2018.03.061.
20. Xu, Q.; Zhang, C.; Wen, C.; Wang, P. A novel composite nonlinear controller for stabilization of constant power load in DC microgrid. *IEEE Trans. Smart Grid* **2017**, doi:10.1109/tsg.2017.2751755.
21. Herrera, L.; Zhang, W.; Wang, J. Stability analysis and controller design of DC microgrids with constant power loads. *IEEE Trans. Smart Grid* **2017**, *8*, 881–888, doi:10.1109/tsg.2015.2457909.

22. Kang, Q.; Feng, S.; Zhou, M.; Ammari, A.C.; Sedraoui, K. Optimal load scheduling of plug-in hybrid electric vehicles via weight-aggregation multi-objective evolutionary algorithms. *IEEE Trans. Intell. Transp. Syst.* **2017**, *18*, 2557–2568, doi:10.1109/tits.2016.2638898.
23. Li, Y.; Li, Y.; Li, G.; Zhao, D.; Chen, C. Two-stage multi-objective OPF for AC/DC grids with VSC-HVDC: Incorporating decisions analysis into optimization process. *Energy* **2018**, *147*, 286–296, doi:10.1016/j.energy.2018.01.036.
24. Adaryani, M.R.; Karami, A. Artificial bee colony algorithm for solving multi-objective optimal power flow problem. *Int. J. Electr. Power Energy Syst.* **2013**, *53*, 219–230, doi:10.1016/j.ijepes.2013.04.021.
25. Yuan, X.; Zhang, B.; Wang, P.; Liang, J.; Yuan, Y.; Huang, Y.; Lei, X. Multi-objective optimal power flow based on improved strength Pareto evolutionary algorithm. *Energy* **2017**, *122*, 70–82, doi:10.1016/j.energy.2017.01.071.
26. Abaci, K.; Yamacli, V. Differential search algorithm for solving multi-objective optimal power flow problem. *Int. J. Electr. Power Energy Syst.* **2016**, *79*, 1–10, doi:10.1016/j.ijepes.2015.12.021.
27. Li, Y.; Wang, J.; Zhao, D.; Li, G.; Chen, C. A two-stage approach for combined heat and power economic emission dispatch: Combining multi-objective optimization with integrated decision making. *Energy* **2018**, *162*, 237–254, doi:10.1016/j.energy.2018.07.200.
28. Shaheen, A.M.; El-Sehiemy, R.A.; Farrag, S.M. Solving multi-objective optimal power flow problem via forced initialised differential evolution algorithm. *IET Gen. Transm. Dis.* **2016**, *10*, 1634–1647, doi:10.1049/iet-gtd.2015.0892.
29. Li, Y.; Feng, B.; Li, G.; Qi, J.; Zhao, D.; Mu, Y. Optimal distributed generation planning in active distribution networks considering integration of energy storage. *Appl. Energy* **2018**, *210*, 1073–1081, doi:10.1016/j.apenergy.2017.08.008.
30. Antonio, L.M.; Coello, C.A.C. Coevolutionary multi-objective evolutionary algorithms: A survey of the state-of-the-art. *IEEE Trans. Evol. Comput.* **2017**, *22*, 851–865, doi:10.1109/tevc.2017.2767023.
31. Sarro, F.; Ferrucci, F.; Harman, M.; Manna, A.; Ren, J. Adaptive multi-objective evolutionary algorithms for overtime planning in software projects. *IEEE Trans. Softw. Eng.* **2017**, *43*, 898–917, doi:10.1109/tse.2017.2650914.
32. Li, W.; Özcan, E.; John, R. Multi-objective evolutionary algorithms and hyper-heuristics for wind farm layout optimisation. *Renew. Energy* **2017**, *105*, 473–482, doi:10.1016/j.renene.2016.12.022.
33. Zhang, X.; Tian, Y.; Jin, Y. A knee point-driven evolutionary algorithm for many-objective optimization. *IEEE Trans. Evol. Comput.* **2015**, *19*, 761–776, doi:10.1109/tevc.2014.2378512.
34. Zhang, X.; Tian, Y.; Cheng, R.; Jin, Y. A decision variable clustering-based evolutionary algorithm for large-scale many-objective optimization. *IEEE Trans. Evol. Comput.* **2018**, *22*, 97–112, doi:10.1109/tevc.2016.2600642.
35. Liu, H.L.; Chen, L.; Zhang, Q.; Deb, K. Adaptively allocating search effort in challenging many-objective optimization problems. *IEEE Trans. Evol. Comput.* **2018**, *22*, 433–448, doi:10.1109/tevc.2017.2725902.
36. Cheng, R.; Jin, Y.; Olhofer, M. Test problems for large-scale multiobjective and many-objective optimization. *IEEE Trans. Cybern.* **2017**, *47*, 4108–4121, doi:10.1109/tcyb.2016.2600577.
37. Cheng, R.; Jin, Y.; Olhofer, M.; Sendhoff, B. A reference vector guided evolutionary algorithm for many-objective optimization. *IEEE Trans. Evol. Comput.* **2016**, *20*, 773–791, doi:10.1109/tevc.2016.2519378.
38. Deb, K.; Jain, H. An evolutionary many-objective optimization algorithm using reference-point-based nondominated sorting approach, part I: Solving problems with box constraints. *IEEE Trans. Evol. Comput.* **2014**, *18*, 577–601, doi:10.1109/tevc.2013.2281535.
39. Ye, X.; Liu, S.; Yin, Y.; Jin, Y. User-oriented many-objective cloud workflow scheduling based on an improved knee point driven evolutionary algorithm. *Knowl. Based Syst.* **2017**, *135*, 113–124, doi:10.1016/j.knosys.2017.08.006.
40. Narimani, M.R.; Azizipanah-Abarghooee, R.; Zoghdar-Moghadam-Shahrekohne, B.; Gholami, K. A novel approach to multi-objective optimal power flow by a new hybrid optimization algorithm considering generator constraints and multi-fuel type. *Energy* **2013**, *49*, 119–136, doi:10.1016/j.energy.2012.09.031.
41. Ye, C.J.; Huang, M.X. Multi-objective optimal power flow considering transient stability based on parallel NSGA-II. *IEEE Trans. Power Syst.* **2015**, *30*, 857–866, doi:10.1109/tpwrs.2014.2339352.
42. Deb, K.; Pratap, A.; Agarwal, S.; Meyarivan, T.A.M.T. A fast and elitist multiobjective genetic algorithm: NSGA-II. *IEEE Trans. Evol. Comput.* **2002**, *6*, 182–197, doi:10.1109/4235.996017.

43. Panda, S.; Yegireddy, N.K. Automatic generation control of multi-area power system using multi-objective non-dominated sorting genetic algorithm-II. *Int. J. Electr. Power Energy Syst.* **2013**, *53*, 54–63, doi:10.1016/j.ijepes.2013.04.003.
44. Zhang, J.; Tang, Q.; Li, P.; Deng, D.; Chen, Y. A modified MOEA/D approach to the solution of multi-objective optimal power flow problem. *Appl. Soft Comput.* **2016**, *47*, 494–514, doi:10.1016/j.asoc.2016.06.022.
45. Yokoyama, R.; Bae, S.H.; Morita, T.; Sasaki, H. Multiobjective optimal generation dispatch based on probability security criteria. *IEEE Trans. Power Syst.* **1988**, *3*, 317–324, doi:10.1109/59.43217.
46. Gaing, Z.L. Particle swarm optimization to solving the economic dispatch considering the generator constraints. *IEEE Trans. Power Syst.* **2003**, *18*, 1187–1195, doi:10.1109/tpwrs.2003.814889.
47. Ghasemi, M.; Ghavidel, S.; Ghanbarian, M.M.; Gharibzadeh, M.; Vahed, A.A. Multi-objective optimal power flow considering the cost, emission, voltage deviation and power losses using multi-objective modified imperialist competitive algorithm. *Energy* **2014**, *78*, 276–289, doi:10.1016/j.energy.2014.10.007.
48. Montoya, F.G.; Baños, R.; Gil, C.; Espín, A.; Alcayde, A.; Gómez, J. Minimization of voltage deviation and power losses in power networks using Pareto optimization methods. *Eng. Appl. Artif. Intell.* **2010**, *23*, 695–703, doi:10.1016/j.engappai.2010.01.011.
49. Bouchekara, H.R.E.H.; Chaib, A.E.; Abido, M.A.; El-Sehiemy, R.A. Optimal power flow using an Improved Colliding Bodies Optimization algorithm. *Appl. Soft Comput.* **2016**, *42*, 119–131, doi:10.1016/j.asoc.2016.01.041.
50. Mohamed, A.A.A.; Mohamed, Y.S.; El-Gaafary, A.A.; Hemeida, A.M. Optimal power flow using moth swarm algorithm. *Electr. Power Syst. Res.* **2017**, *142*, 190–206, doi:10.1016/j.epsr.2016.09.025.
51. Davoodi, E.; Babaei, E.; Mohammadi-ivatloo, B. An efficient covexified SDP model for multi-objective optimal power flow. *Int. J. Electr. Power Energy Syst.* **2018**, *102*, 254–264, doi:10.1016/j.ijepes.2018.04.034.
52. Biswas, P.P.; Suganthan, P.N.; Amaratunga, G.A. Optimal power flow solutions incorporating stochastic wind and solar power. *Energy Convers. Manag.* **2017**, *148*, 1194–1207, doi:10.1016/j.enconman.2017.06.071.
53. Thiele, L.; Miettinen, K.; Korhonen, P.J.; Molina, J. A preference-based evolutionary algorithm for multi-objective optimization. *Evol. Comput.* **2009**, *17*, 411–436, doi:10.1162/evco.2009.17.3.411.
54. Pal, N.R.; Bezdek, J.C. On cluster validity for the fuzzy c-means model. *IEEE Trans. Fuzzy Syst.* **1995**, *3*, 370–379, doi:10.1109/tfuzz.1997.554463.
55. Li, Y.; Li, Y.; Li, G. A two-stage multi-objective optimal power flow algorithm for hybrid AC/DC grids with VSC-HVDC. In Proceedings of the 2017 IEEE Power & Energy Society General Meeting, Chicago, IL, USA, 16–20 July 2017; pp. 1–5.
56. Chang, X.; Wang, Q.; Liu, Y.; Wang, Y. Sparse regularization in fuzzy *c*-means for high-dimensional data clustering. *IEEE Trans. Cybern.* **2017**, *47*, 2616–2627, doi:10.1109/tcyb.2016.2627686.
57. Zang, Y.; Sun, W.; Han, S. Grey relational projection method for multiple attribute decision making with interval-valued dual hesitant fuzzy information. *J. Intell. Fuzzy Syst.* **2017**, *33*, 1053–1066, doi:10.3233/jifs-162422.
58. Li, Y.; Li, Y.; Sun, Y. Online static security assessment of power systems based on lasso algorithm. *Appl. Sci.* **2018**, *8*, 1442, doi:10.3390/app8091442.
59. Li, Y.; Yang, Z. Application of EOS-ELM with binary Jaya-based feature selection to real-time transient stability assessment using PMU data. *IEEE Access* **2017**, *5*, 23092–23101, doi:10.1109/ACCESS.2017.2765626.
60. Yang, Z.; Zhong, H.; Bose, A.; Zheng, T.; Xia, Q.; Kang, C. A linearized OPF model with reactive power and voltage magnitude: A pathway to improve the Mw-only DC OPF. *IEEE Trans. Power Syst.* **2018**, *33*, 1734–1745, doi:10.1109/tpwrs.2017.2718551.
61. Levron, Y.; Guerrero, J.M.; Beck, Y. Optimal power flow in microgrids with energy storage. *IEEE Trans. Power Syst.* **2013**, *28*, 3226–3234, doi:10.1109/TPWRS.2013.2245925.
62. Reddy, S.S. Optimal power flow with renewable energy resources including storage. *Electr. Eng.* **2017**, *99*, 685–695, doi:10.1007/s00202-016-0402-5.